\documentclass[12pt]{article}
\newtheorem{theorem}{Theorem}

\newtheorem{corollary}[theorem]{Corollary}

\newtheorem{example}[theorem]{Example}

\newtheorem{lemma}[theorem]{Lemma}

\newtheorem{proposition}[theorem]{Proposition}
\newtheorem{remark}[theorem]{Remark}

\newenvironment{proof}[1][Proof]{\textbf{#1.} }{\ \rule{0.5em}{0.5em}}
\newdimen\dummy
\dummy=\oddsidemargin
\begin{document}

\date{}
\title{The Hawaiian earring group is topologically incomplete}
\author{Paul Fabel \\
%EndAName
Department of Mathematics \& Statistics\\
Mississippi State University}
\maketitle

\begin{abstract}
If an inverse limit space $X$ is constructed in similar fashion to the
Hawaiian earring, then the fundamental group $\pi _{1}(X,p)$ becomes a
topological group with either of two natural but distinct topologies. In
either case, despite being uncountable, we show $\pi _{1}(X,p)$ is not a
Baire space. The Hawaiian earring in particular provides a premier example
of a Peano continuum whose topological fundamental group, despite regularity
and uncountability, admits no complete metric compatible with its topology.
\end{abstract}

\section{Introduction}

Endowed with the quotient topology from the space of based loops in $X$, the%
\textit{\ }familiar fundamental group\textit{\ }$\pi _{1}(X,p)$\ of a space $%
X$ becomes a \textit{topological }group\textit{. }As a subject in its own
right, the theory of \textit{topological fundamental groups} is in the early
stages of development, and the extent to which the theory will prove useful
remains to be seen. For example, the theory has no extra value in the
context of locally contractible spaces, since if $X$ is locally
contractible, then $\pi _{1}(X,p)$ inherits the discrete topology (\cite{fab}%
) and encodes no more information then algebraic data.

On the other hand, if we also consider spaces that are not locally
contractible then the theory begins to earn its keep. For example both the
algebra \textit{and }topology of $\pi _{1}(X,p)$ is invariant under the
homotopy type of the underlying space (Corollary 3.4 \cite{Biss}).
Consequently, by analyzing the \textit{topology} of $\pi _{1}(X,p),$ the
theory has the capacity to distinguish spaces of distinct homotopy types
when the familiar algebraic homotopy groups fail to do so (see \cite{Fab3}
for examples.)

In any mathematical theory (e.g. the theory of groups, the theory of $n$
manifolds, the theory of Banach spaces, etc...) the investigation of
reasonable conjectures tends to work in tandem with the cataloging of
properties of the fundamental examples. This paper was written with the
latter goal in mind. In particular we investigate the topological
fundamental group of one of the simplest nonlocally contractible spaces, the
Hawaiian earring $HE$. The space $HE$ is the union of a null sequence of
simple closed curves joined at a common point $p$. Hence $HE$ is a one
dimensional metric Peano continuum, and the deletion of one point leaves a
(disconnected) one manifold.

The algebra of $\pi _{1}(HE,p)$ has been studied from various perspectives
by other authors (\cite{Bog2}, \cite{can}, \cite{Can2}, \cite{can3}, \cite
{desmit}, \cite{eda1}, \cite{eda3},\cite{eda},\cite{Fisher}, \cite{mor},\cite
{zas2},\cite{zas3}.) However, relatively little is known about the topology
of $\pi _{1}(HE,p).$

The following tool is indispensable for the investigation of $\pi _{1}(HE,p)$%
. Consider $HE=$ $\lim_{\leftarrow }X_{n}$ as the inverse limit of a nested
sequence of retracts $X_{1}\subset X_{2},...$ where $X_{n}$ is the union of $%
n$ simple closed curves joined at a common point $p,$ and consider the
canonical homomorphism $\phi :\pi _{1}(HE,p)\rightarrow \lim_{\leftarrow
}\pi _{1}(X_{n},p).$ Since $\phi $ is continuous (Proposition 3.3 \cite{Biss}%
), and since $X_{n}$ is locally contractible, the question of whether $\pi
_{1}(HE,p)$ is a $T_{1}$ space (a space whose one point subsets are closed)
is precisely equivalent (see \cite{fab5}) to the question of whether $\phi $
is one to one. A comprehensive affirmative answer to the latter question
constitutes the lion's share of a paper of Morgan and Morrison \cite{mor}.

Since the $T_{1}$ space $\pi _{1}(HE,p)$ is also a topological group, it
follows from general properties of topological groups that $\pi _{1}(HE,p)$
is in fact completely regular. It is an open question whether $\pi
_{1}(HE,p) $ is normal or whether $\pi _{1}(HE,p)$ is metrizable.

Having said this we should point out that the metrizabilty (and hence
normality) of $\pi _{1}(HE,p)$ would follow from the published work of
another author, work that this author believes to be incorrect. A lack of
clarity makes it difficult to pinpoint the location of the first mistake.
The paper \cite{Biss} seems to claim that $\phi :\pi _{1}(HE,p)\rightarrow
\lim_{\leftarrow }\pi _{1}(X_{n},p)$ is an isomorphism with continuous
inverse (Theorem 8.1 \cite{Biss}). Surjectivity is easily seen to fail
since, in $\lim_{\leftarrow }\pi _{1}(X_{n},p),$ the ``infinite word'' $%
x_{1},x_{2},x_{1}^{-1},x_{2}^{-1},x_{1},x_{3},x_{1}^{-1},x_{3}^{-1},...$ has
no counterpart in $\pi _{1}(HE,p).$ The failure of $\phi $ to be a
homeomorphism onto its image is not so easy to see, and an explicit
counterexample is contained in my paper \cite{Fabel}.

Fortunately the main results of this paper do not depend on these
contestable claims. We prove that if $X$ is the Hawaiian earring, or if $X$
is constructed in similar fashion to the Hawaiian earring, then $\pi
_{1}(X,p),$ despite being uncountable, is not a Baire space.

An immediate application is the following: Ignoring the question of whether
the regular space $\pi _{1}(X,p)$ is metrizable, we can conclude by the
Baire category theorem that $\pi _{1}(X,p)$ admits no \textit{complete}
metric compatible with its topology. This in turn yields a `no retraction
theorem' whose roots depend on the \textit{topological} rather than the 
\textit{algebraic }properties of the groups in question: Because $\pi
_{1}(X,p)$ is not a Baire space, the main result of \cite{fab2} shows that $%
X $ cannot be embedded as retract of any space $Y$ such that $\pi _{1}(Y,p)$
admits a complete metric.

Finally, in the process of proving $\pi _{1}(X,p)$ is not a Baire space we
also prove that $\iota m(\phi )$ is not a Baire space in $\lim_{\leftarrow
}\pi _{1}(X_{n},p).$ Thus this paper seems to contain the premier exhibition
of the following phenomenon: There exists a Peano continuum $X$ such that
with either of two natural but distinct $T_{1}$ topological group
structures, the fundamental group of $X$ admits no complete metric
compatible with its topology.

It is hoped that uncovering the properties of the simplest examples such as $%
HE$ will serve as a useful benchmark when pursuing further developments in
the theory of topological fundamental groups.

\section{Preliminaries}

Given a set $S$ a \textbf{word} in $S$ is either the empty set or a finite
sequence in $S.$ The entries of a word are called \textbf{letters. }A 
\textbf{subword} $v$ of the word $w$ is the word obtained by deleting
finitely may letters from $w.$

Suppose $H_{1},H_{2},..$ is a sequence (either finite or infinite) of
pairwise disjoint groups. Let $W$ denote the set of words in the set $%
H_{1}\cup H_{2}\cup ...$ .

If $a$ is a letter of the word $w\in W$ then $a$ is of \textbf{type }$i$ if $%
a\in H_{i}.$

Define $K:W\rightarrow W$ as follows. Suppose $w$ is a word in $W.$ To
define $K(w)$ first partition $w$ into a minimal number of cells such that
each cell is a subword of $w$ consisting of consecutive letters of $w,$ such
that each letter is of the same type. Delete each cell whose multiply to
make $id_{H_{i}}.$ Replace the remaining cells by the product of its
letters. Call the resulting subword $K(w).$ Notice $K(w)$ has no more
letters than $w,$ and that $K(w)$ and $w$ have the same number of letters if
and only if $K(w)=w.$ Thus the sequence $w,K(w),K^{2}(w),...$ is eventually
constant. If $K(w)=w$ we say $w$ is \textbf{reduced.} Moreover a word $w$ is
reduced if and only if $w$ is empty word or $w$ is a nonempty word in the
set $(H_{1}\cup H_{2}\cup ..)\backslash \{id_{H_{1}},id_{H_{2}},..\}$ such
that no two consecutive letters are of the same type.

Define $RD:W\rightarrow W$ such that $RD(w)=\lim_{n\rightarrow \infty
}K^{n}(w).$ Declare two words $w$ and $v$ in $W$ equivalent if $RD(v)=RD(w).$
By construction each equivalence class contains a unique reduced word.

The group \textbf{free product }$H_{1}\ast H_{2}...$ consists of the set of
all reduced words in $H_{1}\cup H_{2}....$with multiplication determined by
the rule $(a_{1},a_{2},..a_{n})\ast
(b_{1},b_{2},..b_{m})=RD(a_{1},a_{2},..a_{n},b_{1},..b_{m}).$ Thus $%
\emptyset $ serves as identity$,$ and inverses are determined such that if $%
w=a_{1},a_{2},...a_{n}$ then $w^{-1}=a_{n}^{-1},.,a_{1}^{-1}.$ See \cite
{Hatch} for more details on free products of groups.

A surjective map $q:Y\rightarrow X$ is a \textbf{quotient map} provided $U$
is open in $X$ if and only if $f^{-1}(U)$ is open in $Y.$

Suppose $X$ is a topological space and $p\in X$ . Endowed with the compact
open topology, let $C_{p}(X)=\{f:[0,1]\rightarrow X$ such that $%
f(0)=f(1)=p\}.$ Then the \textbf{topological fundamental group} $\pi
_{1}(X,p)$ is the quotient space of $C_{p}(X)$ obtained by considering the
path components of $C_{p}(X)$ as points. Endowed with the quotient topology, 
$\pi _{1}(X,p)$ is a topological group (Proposition 3.1\cite{Biss}).
Moreover a map $f:X\rightarrow Y$ determines a continuous homomorphism $%
f^{\ast }:\pi _{1}(X,p)\rightarrow \pi _{1}(Y,f(p))$ (Proposition 3.3 \cite
{Biss}).

If $A_{1},A_{2},...$ are topological spaces and $f_{n}:A_{n+1}\rightarrow
A_{n}$ is a continuous surjection then, (endowing $A_{1}\times A_{2}..$ with
the product topology), the \textbf{inverse limit space }$\lim_{\leftarrow
}A_{n}=\{(a_{1},a_{2},...)\in (A_{1}\times A_{2}...)|f_{n}(a_{n+1})=a_{n}\}.$

A topological space $Z$ is a \textbf{Baire} space provided $\cap
_{n=1}^{\infty }U_{n}$ is dense in $Z$ whenever each of $U_{1},U_{2},..$ is
an open dense subspace of $Z.$ It is a classical fact that every complete
metric space is a Baire space (Theorem 7.2 \cite{Munk}).

\section{Definitions}

For the remainder of the paper we make the following assumptions.

All topological fundamental groups are endowed with the quotient topology.

Suppose $X^{\symbol{94}}=\cup _{n=1}^{\infty }X_{n}$ is a metrizable
topological space. Assume $X_{1}\subset X_{2}\subset X_{3}..$ is a nested
sequence of locally simply connected subspaces with $p\in X_{1}.$

Assume $X_{n+1}=X_{n}\cup Y_{n}$ with $X_{n}\cap Y_{n}=p$ and assume $Y_{n}$
is not simply connected.

Let $r_{n}:X_{n+1}\rightarrow X_{n}$ denote the retraction collapsing $%
Y_{n+1}$ to the point $p$.

Assume $X^{\symbol{94}}$ admits a topologically compatible metric such that $%
diam(Y_{n})\rightarrow 0.$ This assumption serves, for example, to
distinguish the (compact) Hawaiian earring from the (noncompact) union of
countably many circles of radius $1$ joined at the common point $p.$

Let $X=\lim_{\leftarrow }X_{n},$ the inverse limit space under the maps $%
r_{n}.$

Let $R_{n}:X\rightarrow X_{n}$ denote the map $R_{n}(x_{1},x_{2},...)=x_{n}$

Let $\phi :\pi _{1}(X,\{p\})\rightarrow \lim_{\leftarrow }\pi _{1}(X_{n},p)$
denote the canonical homomorphism $\phi ([f])=([R_{1}(f)],[R_{2}(f)]...).$

Let $H_{n}=\pi _{1}(Y_{n},p)$.

Let $G_{n}=H_{1}\ast ...\ast H_{n}$.

Let $V$ denote the set of all words in $\cup _{n=1}^{\infty }H_{n}.$ Words
in $V$ are allowed to be not reduced.

Define $\kappa _{m}:V\rightarrow V$ such that $\kappa _{m}(w)=v$ \ if $v$ is
the subword of $w$ obtained by deleting from $w$ all letters of type $i$
where $i\geq m+1.$

Let $\psi _{n}:G_{n+1}\rightarrow G_{n}$ denote the canonical epimorphism $%
\psi _{n}=RD(\kappa _{n}).$ Thus, to compute $\psi _{n}$ of the reduced word 
$w_{n+1},$ first delete from $w_{n+1}$ all letters of type $n+1,$ and then
reduce the resulting subword as much as possible.

Let $G$ denote the subgroup of $\lim_{\leftarrow }G_{n}$consisting of all
sequences $(w_{1},w_{2},...)$ such that for each $m\geq 1$ the following
sequence is eventually constant: $\kappa _{m}(w_{1}),\kappa _{m}(w_{2}),...$

Define $\sigma :G\rightarrow \{1,2,3,...\}$ such that $\sigma (\{w_{n}\})=N$
if $N$ is minimal such that $\kappa _{1}(w_{N}),\kappa _{1}(w_{N+1}),...$ is
constant.

\section{Facts about the preliminary data}

\begin{proposition}
$X^{\symbol{94}}$ is canonically homeomorphic to $X$ under the map $h:\cup
_{n=1}^{\infty }X_{n}\rightarrow X$ satisfying $%
h(x_{n})=(p,p,,,x_{n},x_{n},..)$. Moreover $h_{X_{n}}\hookrightarrow X$ is
an embedding. Henceforth we consider $X^{\symbol{94}}$ and $X$ to be the
same space.
\end{proposition}

\begin{proposition}
\label{vk}The van Kampen theorem determines a canonical isomorphism $%
j_{n}^{\ast }:\pi _{1}(X_{n},p)\rightarrow G_{n}.$ Henceforth we will
consider these groups to be the same. See \cite{Hatch} for details.
\end{proposition}

Most of \cite{mor} is devoted to Theorem 4.1, which, translated into the
notation of this paper, becomes Proposition \ref{metprop}.

\begin{proposition}
\label{metprop}The canonical map $\phi :\pi _{1}(X,p)\rightarrow
\lim_{\leftarrow }G_{n}$ is a monomorphism. Moreover $im(\phi )=G.$
\end{proposition}

The following proposition is perhaps surprising and is the subject of \cite
{Fabel}.

\begin{proposition}
The continuous isomorphism $\phi :\pi _{1}(X,p)\rightarrow G$ is \textbf{not}
a homeomorphism.
\end{proposition}

\section{Main result: $\protect\pi _{1}(X,p)$ is not a Baire space}

\begin{lemma}
\label{b1}For each $N\geq 1,$ the set $\sigma ^{-1}(N)$ is closed in $%
\lim_{\leftarrow }G_{n}.$
\end{lemma}

\begin{proof}
Recall an element $f\in \lim_{\leftarrow }G_{n}$ is a function $%
f:\{1,2,3,...\}\rightarrow \cup _{n=1}^{\infty }G_{n}$ such that $%
w_{n}=f(n)\in G_{n}$ and $w_{n}$ is a reduced word in $H_{1}\ast ...\ast
H_{n}$ and $\psi _{n}(w_{n+1})=w_{n}.$ Suppose the sequence $\{f_{m}\}\in
\sigma ^{-1}(N)$ converges to $f\in \lim_{\leftarrow }G_{n}.$ Suppose $N\leq
i.$ Since $f_{m}$ converges to $f$ pointwise, and since both $G_{i}$ and $%
G_{i+1}$ have the discrete topology, there exists $M$ such that $%
f(i)=f_{M}(i)$ and $f(i+1)=f_{M}(i+1).$ Since $f_{M}\in \sigma ^{-1}(N)$ we
know $\kappa _{1}(f_{M}(i))=\kappa _{1}(f_{M}(i+1)).$ Thus, if $i\geq N$ $%
\kappa _{1}(f(i))=\kappa _{1}(f(i+1)).$ In similar fashion we may choose $K$
such that $f(N-1)=f_{K}(N-1)$ and $f(N)=f_{K}(N).$ It follows that $\kappa
_{1}(f(N-1))\neq \kappa _{1}(f(N)).$ Thus $\sigma (f)=N.$ Hence $\sigma
^{-1}(N)$ is closed in $\lim_{\leftarrow }G_{n}.$
\end{proof}

\begin{lemma}
\label{b2}For each $K\geq 1$ the set $\sigma ^{-1}(K)$ contains no open
subset of $G.$
\end{lemma}

\begin{proof}
Suppose $\sigma (f)=K.$ For each $i\geq 1,$ choose $h_{i}\in H_{i}\backslash
id_{H_{i}}.$ For $n\geq 2$ define $g_{n}\in \lim_{\leftarrow }G_{n}$ such
that if $i\geq n$ then $g_{n}(i)=h_{1},h_{n},h_{1}^{-1},h_{n}^{-1}.$ .

In similar fashion define $g^{n}\in \lim_{\leftarrow }G_{n}$ such that if $%
i\geq n$ then $g^{n}(i)=h_{n},h_{1},h_{n}^{-1},h_{1}^{-1}.$ $d$

Note $g_{n}\rightarrow (\emptyset ,\emptyset ,...)$ and $g^{n}\rightarrow
(\emptyset ,\emptyset ,..).$ Let $f_{n}=g_{n}\ast f$ and $f^{n}=g^{n}\ast f.$
Since $\lim_{\leftarrow }G_{n}$ is a topological group, $f_{n}\rightarrow f$
and $f^{n}\rightarrow f$. Let $\kappa
_{1}(f(K))=x_{k_{1}},x_{k_{2}},..x_{k_{s}}.$ Suppose $n\geq K+1.$ Since $%
\sigma (f)<n$ we conclude

\[
\kappa _{1}(f(n))=x_{k_{1}},x_{k_{2}},..x_{k_{s}}=\kappa _{1}(f(n-1))=\kappa
_{1}(f_{n}(n-1))=\kappa _{1}(f^{n}(n-1)). 
\]

Next we show that at least one of the following inequalities hold: $\kappa
_{1}(f_{n}(n))\neq x_{k_{1}},x_{k_{2}},..x_{k_{s}}$ or $\kappa
_{1}(f^{n}(n))\neq x_{k_{1}},x_{k_{2}},..x_{k_{s}}.$

Let $f(n)=x_{1},x_{2},...,x_{m}\in H_{1}\ast ..\ast H_{n}.$

If $x_{1}\notin H_{n}$ then $%
f_{n}(n)=h_{1},h_{n},h_{1}^{-1},h_{n}^{-1},x_{1},x_{2}...x_{m}$ and $\kappa
_{1}(f_{n}(n))=h_{1},h_{1}^{-1},x_{k_{1}},x_{k_{2}},..x_{k_{s}}\neq
x_{k_{1}},x_{k_{2}},..x_{k_{s}}$

If $x_{1}\in H_{n}$ then $%
f^{n}(n)=h_{n},h_{1},h_{n}^{-1},h_{1}^{-1},x_{1},x_{2}...x_{m}$ and $\kappa
_{1}(f^{n}(n))=h_{1},h_{1}^{-1},x_{k_{1}},x_{k_{2}},..x_{k_{s}}\neq
x_{k_{1}},x_{k_{2}},..x_{k_{s}}.$

Thus, for each $n$ at least one of $f_{n}$ or $f^{n}$ does not belong to $%
\sigma ^{-1}(K).$ Hence $\sigma ^{-1}(K)$ contains no open subset of $G.$
\end{proof}

\begin{example}
\label{trivex}Replace the usual topology of the rational numbers $Q$ with
the discrete topology creating the space $Q^{\symbol{94}}.$ Then $id:Q^{%
\symbol{94}}\rightarrow Q$ is continuous and one to one but not an
embedding. In particular the set $\{0\}\subset Q$ contains no open set.
However $\{0\}$ is open in $Q^{\symbol{94}}.$
\end{example}

\begin{remark}
It is shown in \cite{Fabel} that $\phi $ is not an embedding. Hence,
considering example \ref{trivex}, Lemma \ref{b3} is not an immediate
consequence of Lemma \ref{b2}. However the same idea drives both proofs.
\end{remark}

\begin{lemma}
\label{b3}For each $K\geq 1$ the set $\phi ^{-1}((\sigma ^{-1}(K))$ contains
no open subset of $\pi _{1}(X,\{p\}).$
\end{lemma}

\begin{proof}
Since $Y_{n}$ is not simply connected for each $n\geq 1$ select $h_{n}\in
C_{p}(X_{n})$ such that $[h_{n}]\in \pi _{1}(Y_{n,}p)\backslash id_{n}.$ Let 
$h_{\infty }$ denote the constant map $p.$ Let $U=\phi ^{-1}(\sigma
^{-1}(K)).$ Let $F:C_{p}(X)\rightarrow \pi _{1}(X,p)$ denote the canonical
surjective quotient map. Let $A\subset U.$ Let $B=F^{-1}(U).$ Since $F$ is a
quotient map, in order to prove $A$ is not open it suffices to prove $B$ is
not open in $C_{p}(X).$ Suppose $[\alpha ]\in A.$ For each $n\in \{\infty
\}\cup \{1,2,3,...\}$ let $g_{n}=h_{1}\ast h_{n}\ast h_{1}^{-1}\ast
h_{n}^{-1}\ast \alpha $ denote the path in $C_{p}(X)$ obtained by
concatenation and compatible with the partition $\{0,\frac{1}{5},\frac{2}{5},%
\frac{3}{5},\frac{4}{5},1\}.$

In similar fashion define $g^{n}=h_{n}\ast h_{1}\ast h_{n}^{-1}\ast
h_{1}^{-1}\ast \alpha .$ Note $g_{n}\rightarrow g_{\infty }$ and $%
g^{n}\rightarrow g^{\infty }$ uniformly and $\{g^{\infty },g_{\infty
}\}\subset B$ since $g_{\infty }$ and $g^{\infty }$ are path homotopic in $X$
to $\alpha .$ For $n\geq 2$ let $f^{n}=$ $\phi F(g^{n})$, let $f_{n}=\phi
F(g_{n})$ and let $f=\phi F(\alpha ).$ Since $\phi F$ is continuous, we
conclude $f_{n}\rightarrow f$ and $f^{n}\rightarrow f.$ Now we may appeal to
the proof of Lemma \ref{b2} to conclude that if $n\geq K+1$ then at least
one $f^{n}$ or $f_{n}$ is not an element of $\sigma ^{-1}(K).$ Thus at least
one of $g^{n}$ or $g_{n}$ does not belong to $B.$ Hence $B$ is not open and
this proves Lemma \ref{b3}.
\end{proof}

\begin{theorem}
The topological groups $G$ and $\pi _{1}(X,p)$ are not Baire spaces.
\end{theorem}

\begin{proof}
To prove that a topological space $Z$ is note a Baire space it suffices to
prove $Z$ is the countable union of closed subspaces $Z_{1}\cup Z_{2}...$%
such that $Z_{n}$ has empty interior for each $n.$ Taking $Z_{n}=\sigma
^{-1}(n)$ it follows directly from Lemmas \ref{b1} and \ref{b2} that $G$ is
not a Baire space. To prove $\pi _{1}(X,p)$ is not a Baire space consider
the sets $A_{n}=\phi ^{-1}(Z_{n}).$ Since $\phi $ is continuous, $A_{n}$ is
closed by Lemma \ref{b1}. Lemma \ref{b3} shows $A_{n}$ has empty interior.
\end{proof}

\begin{corollary}
Endowed with either the quotient topology or the inverse limit space
topology, $\pi _{1}(X,p)$ does not admit a complete metric.
\end{corollary}

\end{document}